\documentclass{amsart}
\usepackage{amsfonts}
\usepackage{enumerate}
\usepackage{graphicx}

\newtheorem{theorem}[subsection]{Theorem}

\newtheorem{definition}[subsubsection]{Definition}

\numberwithin{equation}{subsection}

\title{An awkward graph}
\author{Mary Rees}
\begin{document}
\maketitle

\begin{abstract}
Given a rational map $f:\overline{\mathbb C}\to \overline{\mathbb C}$ and a finite graph $G\subset \overline{\mathbb C}$ such that $f(G)\subset G$ and $f$ is expanding on some neighbourhood of $G$, we show that there is another finite graph $G'\subset \bigcup _{n\ge 0}f^{-n}(G)$ in an arbitrarily small neighbourhood of $G$ such that $f^N(G')\subset G'$ for some integer $N$ but $\bigcup _{i=0}^{N-1}f^{i}(G')$ contains accumulating {\em{plaits}} and {\em{nests}}
\end{abstract}

\section{}\label{1}
Let $X$ be a metric space with metric $d$. A map $f:X\to X$ is {\em{expanding}} if there are $\lambda >1$ and an integer $N$ such that $d(f^N(x),f^N(y))\ge \lambda d(x,y)$ for all $x$, $y\in X$. The construction of  Markov partitions for a continuous expanding map of a compact metric space is one of the most basic in dynamical systems. The construction is given, for example, in \cite{P-U}. If $X$ is a compact surface, the existence of a Markov partition for $X$ is closely related to existence of a closed invariant set $Y\subset X$ with empty interior and satisfying $f(Y)\subset Y$. The extent to which the existence of such a set $Y$ is related to the existence of a Markov partition for $Y$ depends on the precise definition of Markov partition. But clearly it would be desirable to have $Y$ being a finite graph, for example: a property which has to be worked for. This paper grew out of work to construct an invariant graph for a  hyperbolic rational map. The following theorem was proved by  F.T. Farrell and L.E. Jones in the 1970's, in slightly different language.

\begin{theorem}\cite{F-J}\label{1.1} Let $f:S\to S$ be an expanding map of a compact surface $S$. Let $G_0\subset S$ be any finite graph subject to mild conditions (trivalent, no free vertices, closures of complementary components are closed topological discs, the boundaries of any two complementary components intersect in at most one edge) and let $\varepsilon >0$ be given. Then there exists a graph $G$ which is isotopic to $G_0$ via an isotopy within $\varepsilon $ of the identity, and an integer $N>0$ such that $f^N(G)\subset G$.\end{theorem}

An interesting special case of this result is studied in \cite{B-M}. There it is proved (Theorem 1.2) that if $f:S^2\to S^2$ is an {\em{expanding Thurston map}} then there exists a Jordan   curve ${\mathcal{C}}$ and an integer $N$ such that $f^{-N}({\mathcal{C}})={\mathcal{C}}$. In fact, ${\mathcal{C}}$ can be chosen to be isotopic, relative to the postcritical set, to any Jordan curve passing through the postcritical set -- assuming that $N$ is large enough. A number of examples are given. In particular, Example 15.17 of \cite{B-M}  illustrates the care needed to ensure that ${\mathcal{C}}$, which is constructed via an iterative process, does not have self-intersections.

Suppose that $f$ is an expanding map of a compact surface and $G$ is a finite invariant graph with $f^N(G)\subset G$ for some integer $N$. Then $G\subset f^{-N}(G)$ and if $G^0=\bigcup _{i=0}^{N-1}f^{-i}(G)$ then $G^0\subset f^{-1}(G^0)$, but $G^0$ might not be a finite graph. A map on a surface is very frequently expanding only on some open subset   $U$ of the surface such that $\overline{U}\subset f(U)$. Theorem \ref{1.1} adapts to this situation. One such adaptation applies to hyperbolic rational maps and to rational maps for which every critical point is attracted to an attractive or parabolic  periodic orbit, maps which are sometimes called geometrically finite. In particular the following result was proved in \cite{R1} (Corollary 1.2 of that paper). 

\begin{theorem}\label{1.2} Let $f:\overline{\mathbb C}\to \overline{\mathbb C}$ be a rational map with connected Julia set $J$, such that the forward orbit of each critical point is attracted to an attractive or parabolic periodic orbit, and such that the closure of any Fatou component is a closed topological disc, and all of these are disjoint. Then there exists a  connected graph $G'\subset \overline{\mathbb C}$ such that the following hold.
\begin{enumerate}
\item  $G'\subset f^{-N}(G')$. 
\item $G'$ does not intersect the closure  of any periodic Fatou component.
\item Any component of $\overline{\mathbb C}\setminus G'$ contains at most one periodic Fatou  component of $f$.
  \end{enumerate}
  \end{theorem}
  
  Actually, the result in \cite{R1} is that we can take $N=1$.  That is more difficult. Once again, if $G'\subset f^{-N}(G')$, there is no guarantee that $G^0=\bigcup _{i=0}^{N-1}f^{-i}(G')$ is a finite graph, although it does satisfy $G^0\subset f^{-1}(G^0)$.  It was shown in \cite{R1} that there is, however, a finite graph $G\subset \overline{\bigcup _{n\ge 0}f^{-n}(G')}$ in an arbitrarily small neighbourhood of $G^0$. In the process, it was realised that there are two quite separate mechanisms which cause awkward intersections between the different  finite graphs $f^{-i}(G')$  for $0\le i<N$. In this paper these mechanisms will be called {\em{plaiting}} and {\em{nesting}}
  
  \subsection{Preliminary examples of plaiting and nesting}\label{1.3}
  
  We will start by looking at plaiting and nesting for finitely many semi-infinite arcs in $\mathbb C$, each homeomorphic to $[0,\infty )$, and intersecting the common endpoint, carried by homeomorphism to $0$. Also, in both cases, we first consider the lifts of the structure under the exponential map.
  
 Fix an integer $N$ and $a\in \mathbb R$ and for $0\le k<N$ and $n\in \mathbb Z$ let $\Gamma _{k,n}$ be the arc
  $$\{ (x+ai\sin (x-2\pi k/N)+\pi in):x\in \mathbb R\} .$$
 Then
 $$\Gamma _{k,n}+\frac{2\pi }{N}=\begin{array}{ll}\Gamma _{k+1,n}&\mbox{ if }k<N-1\\ \Gamma _{0,n}&\mbox{ if }k=N-1\end{array}$$ 
 Then
 $$\gamma _k=\exp (\Gamma _{k,0}) \cup \{0\} =\exp(\Gamma _{k,2n})\cup \{ 0\} $$ is an infinite arc, homeomorphic to $[0,\infty )$, for any integer $n$, as is
 $$-\gamma _k=\exp (\Gamma _{k,2n+1})$$
 for any integer $n$. We have
 $$e^{2\pi /N}\cdot \gamma _k=\begin{array}{ll}\gamma _{k+1}&\mbox{ if }k<N-1\\ \gamma _0&\mbox{ if }k=N-1\end{array}$$

 All the arcs $\gamma _k$ intersect at $0$, and for all values of $a\ne 0$, for any $k\ne \ell $, there are also infinitely many other intersection points accumulating at $0$.  The way in which these intersection points occur gives rise to either {\em{plaiting}} or {\em{nesting}}.

$\Gamma _{k,n}\cap \Gamma _{\ell ,m}\ne \emptyset $ if, for some $x$,
$$a\sin (x-2\pi k/N)+\pi n=a\sin (x-2\pi \ell/N)+\pi m,$$
that is if 
\begin{equation}\label{1.3.1}a(\sin (x-2\pi k/N)-\sin (x-2\pi \ell /N))=\pi (m-n).\end{equation}
The largest modulus value of 
$\sin (x-2\pi k/N)-\sin (x-2\pi \ell /N)$ occurs when 
$$\cos (x-2\pi k/N)=\cos (x-2\pi \ell /N),$$
that is, when
$$\sin \left( x-\pi \frac{k+\ell}{N}\right) \sin \left( \pi \frac{k-\ell}{N}\right) =0,$$
that is, when $x=\pm \pi (k+\ell ) /N\mbox{ mod }2\pi $. 

So
$$\mbox{Max}_x\left\vert a( \sin (x-2\pi k/N)-\sin (x-2\pi \ell /N))\right\vert = 2\left\vert a\sin \left( \frac{\pi (k-\ell)}{N}\right) \right\vert $$
So 
$$\mbox{Max}_{x,k,\ell}\left\vert a( \sin (x-2\pi k/N)-\sin (x-2\pi \ell /N))\right\vert =\begin{array}{lll}2|a|&\mbox{ if }&N\mbox{ is even,}\\
2|a|\sin \left( \frac{\pi}{2}-\frac{\pi }{2N}\right) &\mbox{ if }&N\mbox{ is odd.}\end{array}$$
So if $N$ is even and $|a|<\pi /2$, or $N$ is odd and $|a\sin \pi ((N-1)/2N)|<\pi/2 $, then intersections between the arcs  $\gamma _k$ and $\gamma _\ell $ occur in the same order on $\gamma _k$ and $\gamma _\ell $, for all $k$ and $\ell $. This is {\em{plaiting}}, of which the formal definition will be given below. Whatever the value of $a\ne 0$, intersections between $\Gamma _{k,n}$ and $\Gamma _{\ell ,n}$ occur whenever
$$\sin (x-2\pi k/N)=\sin (x-2\pi \ell /N),$$
which happens twice in every closed interval of length $2\pi $. This gives intersections between $\gamma _k$ and $\gamma _\ell $ accumulating on $0$. 

If $|a|\ge \pi /2$ for  $N$ even, or  $|a\sin \pi ((N-1)/2N)|\ge \pi/2 $ for $N$ odd, then $\Gamma _{k,n}$ and $\Gamma _{\ell ,n+1}$ intersect when $k-\ell =\pm N/2$ if $N$ is even, or $k-\ell =(N-1)/2$ if $N$ is odd. Then for such $k$ and $\ell $, $\gamma _k\cup \gamma _\ell $  bounds a disc containing $0$ in its interior. This is {\em{nesting}}, of which the formal definition will be given below. The larger the value of $a$, the more values of $(k,n)$ and $(\ell ,m)$ give solutions to (\ref{1.3.1}) with $n\ne m$. For example, if $|a\sin (2\pi /N)|\ge 1$ then there is a solution for $(k,n)$ and $(\ell ,n)$ whenever $k\ne \ell $. 

We note that the transformation $(x,y)\mapsto (x+2\pi /N,y)$ maps $\Gamma _{k,n}$ to $\Gamma _{k+1,n}$, and therefore the transformation $z\mapsto e^{2\pi /N}z$ maps $\gamma _k$ to $\gamma _{k+1}$. So we have local invariance. 

Now we give the formal definitions of plaiting and nesting. The definitions are for arcs ending at a point.
\begin{definition} Suppose that $\gamma _k:[0,\infty )\to S$, for $0\le k<N$ is a semi-infinite arc on a surface $S$ such that all the $\gamma _k$ have a common endpoint $s\in S$, that is, $\gamma _k$ is a homeomorphism and $\gamma _k(0)=s$ for all $k$  Let $\Gamma _{k,n}=\Gamma _{k,0}+2\pi in$ ($n\in \mathbb Z$) be the lifts of $\gamma _k$ under the exponential map, that is, $\exp (\Gamma _{k,n}=\gamma _k$ for all $n\in \mathbb Z$. Then the $\gamma _k$ ($0\le k<N$) are {\em{plaited}} (near $s$) if $\gamma _k\cap \gamma _\ell $ and $\gamma _k\setminus \gamma _\ell $ both accumulate on  $s$ for all $k\ne \ell $  and, replacing $\Gamma _{\ell ,0}$ by $\Gamma _{\ell ,n_\ell }$ for each $\ell $ if necessary, all intersections lift to $\Gamma _{k,0}\cap \Gamma _{\ell ,0}$, for each $k\ne \ell $. Equivalently, $\bigcup _k\gamma _k$ does not bound any disc containing $s$ in its interior.

The $\gamma _k$ are {\em{nested}} (around $s$) if  $\gamma _k\cap \gamma _\ell $ and $\gamma _ k$ are  not plaited, that is, for at least one $k$ there are $\ell $ and $m$ such that $\Gamma _{k,0}\cap \Gamma _{\ell ,0}\ne \emptyset $ and $\Gamma _{\ell ,0}\cap \Gamma _{m,0}\ne \emptyset $ but $\Gamma _{k,n}\cap \Gamma _{m,0}\ne \emptyset $ for some $n\ne 0$. 
\end{definition}

Note that the $\gamma _k$ in the example given of plaiting are plaited in the common language sense, if we 
choose appropriate over and under-crossings --- which will not be displayed, as it has nothing to do with the matter in hand.  
\begin{center}
\includegraphics[width=11cm]{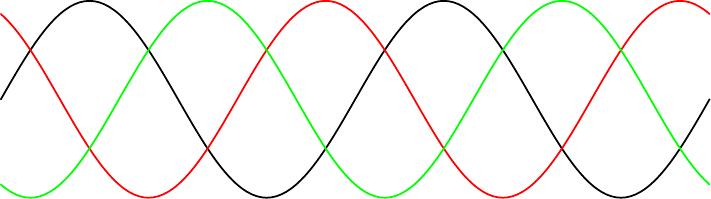}

$\Gamma _{k,0}$ for $N=3$, $0\le k\le 2$ 
\end{center}
\begin{center}
\includegraphics[width=11cm]{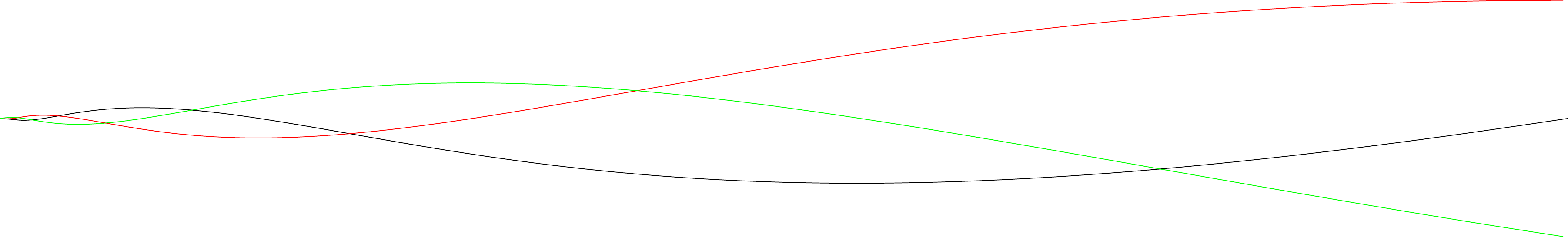}

$\gamma _{k,0}$ for $N=3$, $0\le k\le 2$ and small $a$ (plaiting)
\end{center}

\begin{center}
\includegraphics[width=11cm]{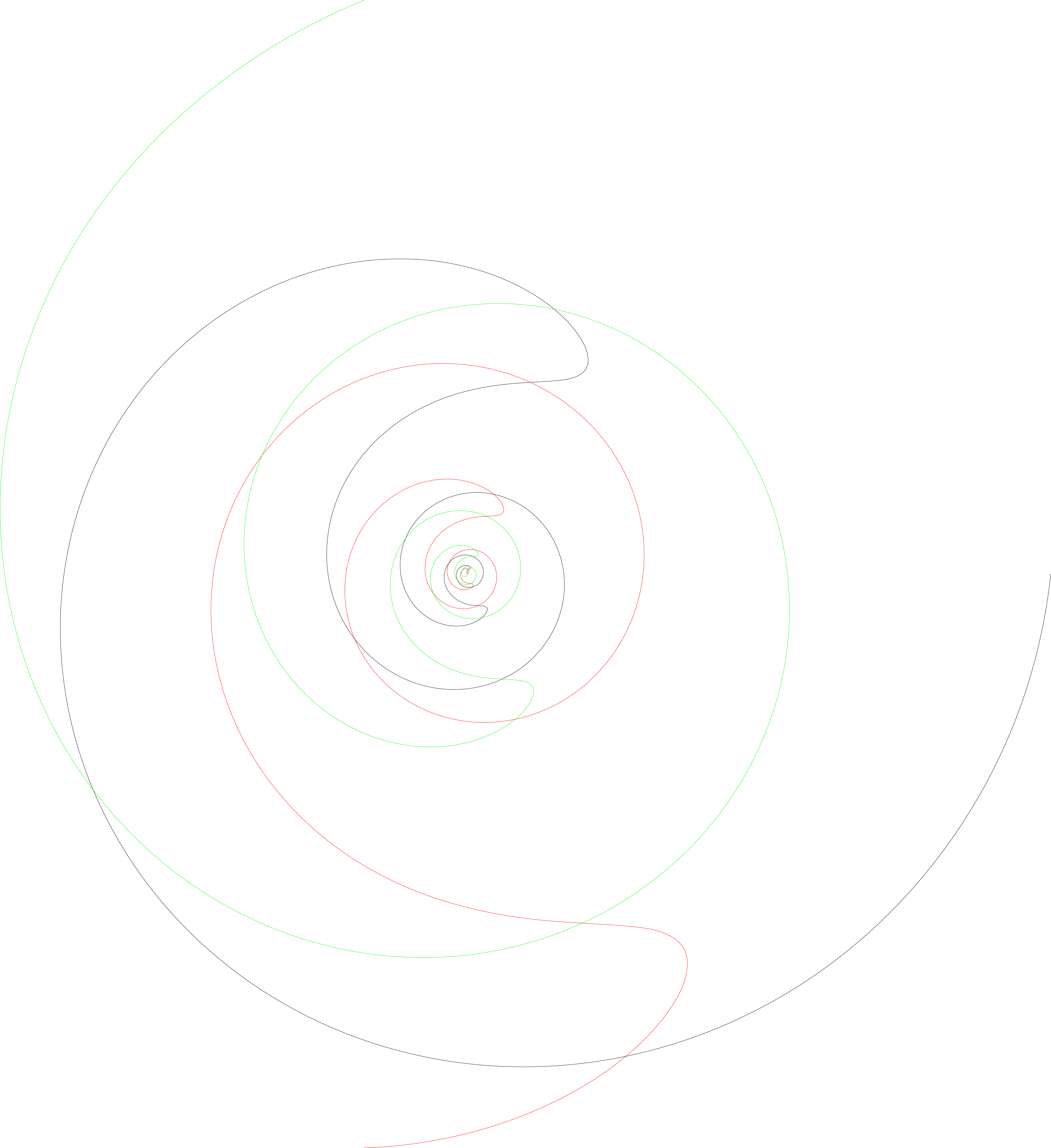}

$\gamma _{k,0}$ for $N=3$, $0\le k\le 2$ and large $a$ (nesting)
\end{center}

Now the theorem that we wish to prove is the following. 
\begin{theorem}\label{1.4} Given a rational map $f:\overline{\mathbb C}\to \overline{\mathbb C}$ and a finite graph $G\subset \mathbb C$ such that $f(G)\subset G$ and $f$ is expanding on some neighbourhood of $G$,  there is another finite graph $G'\subset \bigcup _{n\ge 0}f^{-n}(G)$ in an arbitrarily small neighbourhood of $G$ such that $f^N(G')\subset G'$ for some integer $N$ but $\bigcup _{i=0}^{N-1}f^{i}(G')$ contains  {\em{plaits}}  near a Cantor set, and {\em{nests}} around a Cantor set.\end{theorem}

As is probably already apparent, nesting and plaiting are constructed very similarly, basically by just changing a parameter, although topologically they are different.

\section{The graph examples}\label{2}

\subsection{The basic construction }\label{2.1}
Let $G\subset \overline{\mathbb C}$ be any invariant  connected graph for a rational map $f$, that is, $f(G)\subset G$. Suppose that $f:G\to G$ is expanding in a neighbourhood of $G$, with respect to some metric which is Lipschitz equivalent to the spherical metric. Let $U\subset \overline{\mathbb C}$ and an integer $r$ be such that:
\begin{itemize}
\item $U$ is a neighbourhood of some arc of $G$ such that $f$ is expanding on any component of $f^{-n}(U)$ which intersects $G$, for any $n\ge 0$;
\item $U\subset f^r(U)$, but $U\cap f^i(U)=\emptyset $;
\item  each component of $f^{-r}(U)$ is either contained in $U$ or disjoint from $U$. 
\end{itemize}
Then for a sufficiently large integer $s$ and $N=sr$,  we choose $G'\subset f^{-N}(G)$ so that
$$G'\setminus \bigcup _{k\ge 0}f^{-kr}(U)=G\setminus \bigcup _{k\ge 0}f^{-kN}(U).$$
We define
$$W=U\setminus f^{-r}(U).$$
Then the sets $f^{-kr}(W)$ are all disjoint. If we write $W_{kr}$ for the union of components of $f^{-kr}(W)$ which intersect $G$. Since $f$ is expanding in $U$, all accumulation points of $W_{kr}$ are in $G$. We choose $G'$ so that:
\begin{enumerate}
\item $G'\cap W\subset f^{-N}(G)$;
\item $f^{kr}(W_{kr}\cap G')=W\cap G'$ for $k>0$;
\item if $U_1$ is a component of $f^{-n}(U)$ and $U_1\cap G\ne \emptyset $, then $f^n(U_1\cap G')=U\cap G'$;
\item $G'\cap U$ is an arc.
\end{enumerate}
By 1, 2 and 3,  $f^N(G')\subset G'$ and 
$$G'\subset \bigcup _{k\ge 0}f^{-kN}(G).$$
By 4, $G'$ is a finite graph. Condition 4 can be ensured by choosing $G'\cap W$ so that $G'$ intersects each component of $\partial W$ in exactly two points. Also by 1, 2 and 3, the only transverse intersections between $G'$ and $f^i(G')$ are images under $f^i$, for $0\le i<N$, of intersections in $U$. 
We can choose $G'\cap W_r$ so that $G'$ has transversal intersections with $f^r(G')\cap W$. Similarly, if $s$ is sufficiently large given $t$ we can choose $G'\cap W_{ir}$ suitable for $i\le t$,  so that $G'$ has transversal intersections with $f^{ir}\cap G'$  for $i\le t$. Then $G'\cap U$ has infinitely many intersections with $f^r(G')$ (and more generally with $f^{ir}(G')$ for $i\le t$), which accumulate at points of $\bigcup_{m\ge 0}\bigcap _{k\ge m}W_{kr}$. For different numbers of components of $U\cap f^{-r}(U)$, this set of accumulation points can be a single point or a Cantor set. We can also choose $G'$ so that the graphs $f^i(G')$  are plaited or nested round the points of  $\bigcup_{m\ge 0}\bigcap _{k\ge m}W_{kr}$.

To choose $f^r(G')\cap G'$ to have just finitely many accumulation points, and just one point in $U$, choose $U$ to be a neighbourhood of a point $a\in G$ of period $r$ under $f$. Let $S$ be the local inverse of $f^r$ which fixes $a$.  We can choose $U$ so that $U\cap f^{-r}(U)=S(U)$.  Choose $N$ as large as desired, and divisible by $r$, with $N=sr$.  We can choose $G'\cap U$ so that $G'\cap U$ intersects $f^{ir}(G')\cap U$ transversally only in $\bigcup _{i\ge 0}S^i(U\setminus f^{-r}(U))$, by choosing $G'$ appropriately in any components of $U\cap f^{-ir}(U)\setminus S^i(U)$ for each $i\ge 2$.  By the choice of $G'\cap U\setminus f^{-N}(U)$, we can obtain either plaiting or nesting of $f^{ir}(G')$  round $a$, for $0\le ir<N$. 

We can choose a larger set of accumulation points of plaiting/nesting by choosing $U$ so that $f^{-r}(U)\cap U$ has more than one component for suitable $r$. Suppose that $a_1$ and $a_2$ are periodic points in $U$ of periods $r_1$ and $r_2$ respectively. Let $r$ be the least common multiple of $r_1$ and $r_2$. Suppose $N$ is divisible by $r$ --- and hence by $r_1$ and $r_2$ --- with $N=sr$.
 Let $S_j$ be the local inverse of $f^{r}$ fixing $a_j$. We can choose $G'\cap (U\setminus (S_1(U)\cup S_2(U))$ so that $f^{i}(G')$, for $0\le i<N$ are plaited or nested round
$$\bigcup _{n\ge 1}\{ S_{i_1}\circ \cdots \circ S_{i_n}(U):i_j\in\{ 1,2\} ,\ 1\le j\le n\} .$$
Once again, if $U\cap f^{-N}(U)$ is larger than $U\setminus (S_1^{s}(U)\cup S_2^{s}(U))$ then we can choose $G'\cap U\setminus (S_1^{s}(U)\cup S_2^{s}(U))$ so that transversal intersections between $f^i(G')\cap U\setminus (S_1^{s}(U)\cup S_2^{s}(U))$ and $G'\cap U\setminus (S_1^{s}(U)\cup S_2^{s}(U))$  occur only in $U\setminus f^{-N}(U)$.

\subsection{Example of nesting round a Cantor set}\label{2.2}

We consider $U$ and $G$ such that $f^{-r}(U)\cap U$ has at least two components and $S_j:U\to U$ are distinct local inverses of $f^{r}$ for $j=1$, $2$, and $S_j:U\to U$ is a contraction. Then we perform the construction described in \ref{2.1} with $N=2r$, to obtain nesting of $G'$ and $f^r(G')$ round a Cantor set. In the following pictures the arcs $\gamma _0$ and $\gamma _1$ are $G'\cap U$ and approximations to $f^r(G')\cap U$ respectively, with $\gamma _0$ coloured black and $\gamma _1$ coloured red.  The pictures show $\gamma _0$ and $\gamma _1$ up to homeomorphism, which means that $\gamma _0$ can be, and is,  represented as a straight line. The $n+1$'st approximation $\gamma _{1,n+1}$ to $\gamma _1$ is obtained from the $n$'th by applying $S_1$ and $S_2$ to $\gamma _0\cup \gamma _{1,n}$. So $\gamma _{1,n+1}$ differs from $\gamma _{1,n}$ only in the images of the maps $S_{i_1}\cdots S_{i_{n+1}}$. There are of course $2^{n+1}$ such inverse images, and two within the image of each map $S_{i_1}\cdots S_{i_n}$. The intersections of $\gamma _{1,n}$ with with the image of each $S_{i_1}\cdots S_{i_n}$ are the same up to translation, in the pictures. The change form $\gamma _{1,n}$ to $\gamma _{1,n+1}$  is  determined by the action of $S_{i_1}$ on $\mbox{Im}(S_{i_2}\cdots S_{i_{n+1}}$. As we will see, the change from $\gamma _{1,n}$ to $\gamma _{1,n+1}$ in each of the $2^n$ different pieces is the same, up to homeomorphism preserving black and red -- and even up to similarity, in the pictures --- as the change from $\gamma _{1,n+4}$ to $\gamma _{1,n+5}$, for all $n$. 
In the Stage 1 picture, one of the two rectangles in which the picture changes at Stage 2, is outlined. In each of the other pictures one of the two rectangles, in which the change has been effected from the previous stage, is outlined, and one of the two rectangles within that in which the picture changes, for the next stage, is also outlined.

 In the limit, intersections between $\gamma _0$ and $\gamma _1$ accumulate on a Cantor set. 
\begin{center}
\includegraphics[width=10cm]{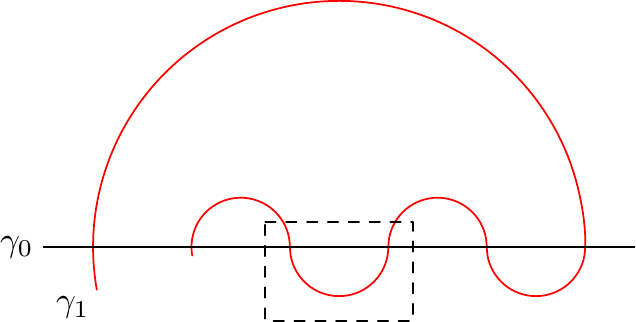}

Stage 1 \end{center}

\vskip 1 true cm

\begin{center}
\includegraphics[width=10cm]{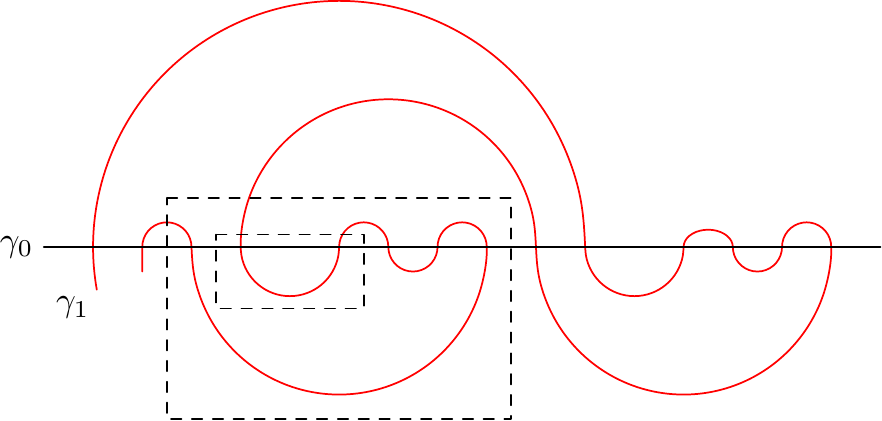}

Stage 2 \end{center}

\vskip 1 true cm

\begin{center}
\includegraphics[width=10cm]{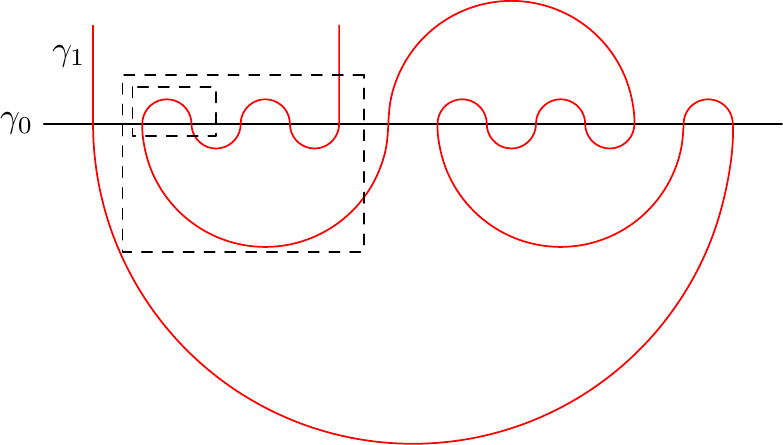}

Stage 3
\end{center}

\vskip 1 true cm

\begin{center}
\includegraphics[width=10cm]{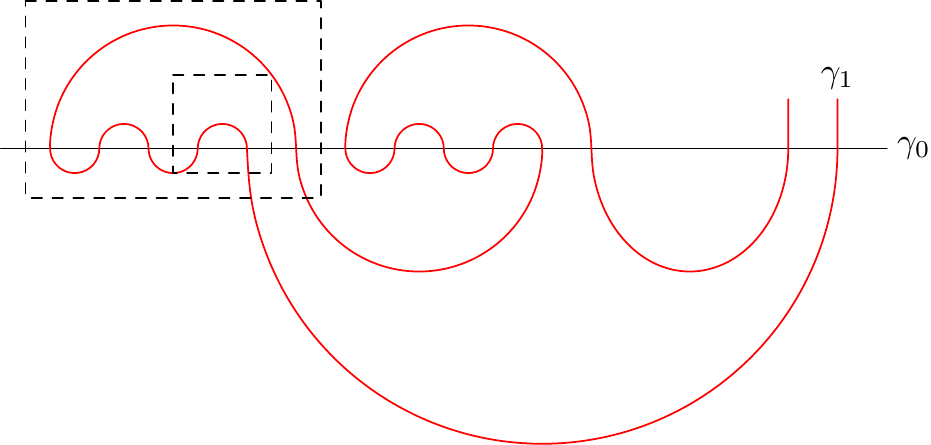}

Stage 4
\end{center}

\end{document}